\documentclass[12pt]{article}
\usepackage{amssymb}

\begin{document}
\begin{center}
{\Large\bf On a generalization 
of the Dvoretzky-Wald-Wolfowitz theorem with an application to a robust optimization problem}  
\end{center}

\begin{center}
{\bf \large  Anna Ja\'skiewicz$^a$, Andrzej S. Nowak$^b$} 
\end{center}{\footnotesize
\noindent$^a$Faculty of Pure and Applied Mathematics,  Wroc{\l}aw University of Science and Technology, 
Wroc{\l}aw, Poland, email: anna.jaskiewicz@pwr.edu.pl\\
\noindent$^b$Faculty of Mathematics, Computer Science and Econometrics, University of Zielona G\'ora,
Zielona G\'ora, Poland, email: a.nowak@wmie.uz.zgora.pl}
\begin{center}
\today
\end{center}
\vspace{3mm}

\noindent{\bf Abstract.} A generalization of the Dvoretzky-Wald-Wolfowitz theorem to the case of 
conditional expectations is provided assuming that the $\sigma$-field on the state space 
has no conditional atoms.\\


\section{Introduction}

Dvoretzky, Wald and Wolfowitz \cite{dww1,dww2} showed that 
for any measurable function taking values in a finite dimensional simplex
there exists a measurable function taking values in the extreme points of the simplex and such that it has 
an integral identical to that of the given function with respect to a finite set of bounded non-atomic measures. 
Its proof is based upon Lyapunov's theorem \cite{l}.
Both results play a crucial role in control theory, statistical decision problems and game theory, see for instance 
\cite{b3,b5,b,fp2,f,krs,ks}. Next, in 1976 Dynkin and Evstigneev \cite{de} provided 
a generalization of Lyapunov's theorem \cite{l}
in the sense that they considered  conditional expectations of set-valued functions assuming 
that the $\sigma$-field on the state space has no conditional atom. 
Recently, this result was applied to stochastic games by He and Sun \cite{hs}. 
They showed that a stochastic game with a so-called coarser transition kernel possesses 
a stationary Markov perfect equilibrium. 

In this paper, we provide a generalization of the Dvoretzky-Wald-Wolfowitz theorem 
to the case with conditional expectations
using the Dynkin and Evstigneev approach.
The proof combines the Carath\'eodory theorem and certain ideas used in \cite{de}.  
Moreover, our result can be viewed as generalizations 
of the Dworetzky-Wald-Wolfowitz reported in \cite{fp2} and \cite{b4}, see Remark 2 for further comments. 
Finally, we illustrate how to use the purification principle 
to solve a robust variational problem with integral constraints.
In contrast to works \cite{art,ap, b1, bl}, we allow the malevolent nature to choose a parameter
whose influences on the distribution of the state of the world.  

\section{Main Purification Result}
Let $(\Omega, \cal F, \mu)$ be a complete probability space with a non-atomic measure $\mu$ and 
let $X$ be a complete separable metric space endowed with the $\sigma$-field ${\cal B }(X)$ of all its Borel subsets. 
By $P(X)$ we denote the space of all  probability measures on $X$ equipped with the topology of weak convergence
and the Borel  $\sigma$-field ${\cal B }(P(X)).$ Let $A$ be a correspondence from $\Omega$ to $X$ such that
$A(\omega)\not= \emptyset$ for each $\omega\in \Omega$ 
and its graph $Gr(A)= \{(\omega,x): \omega\in\Omega, x\in A(\omega)\}$ 
belongs to the product $\sigma$-field ${\cal F}\otimes {\cal B}(X).$ A randomized decision function is a 
measurable mapping $\varphi:\Omega\to P(X)$ such that 
$\varphi(\omega)(A(\omega))=1$ for all $\omega\in\Omega.$ 
The set of all randomized decision functions is denoted by $\Phi.$ 
A non-randomized or pure decision function is
a measurable mapping $f:\Omega\to X$ such that $f(\omega)\in A(\omega)$ for all $\omega\in\Omega.$
Clearly, the set $F$ of all non-randomized decision functions can be viewed as a strict subset of
$\Phi.$ By the von Neumann-Aumann measurable selection theorem 
(see Theorem 18.26 in \cite{ab} or Theorem III.22 in \cite{cv}) 
we know that $F\not=\emptyset.$

Let $\mathbb{R}^n$ be the Euclidean space endowed with the usual norm $\|\cdot\|$ 
and the Borel $\sigma$-field ${\cal B}(\mathbb{R}^n).$ 
Consider a family of functions $g_i:Gr(A)\to \mathbb{R}$ ($i=1,...,n$), 
which are ${\cal F}\otimes {\cal B}(X)$-measurable and 
\begin{equation}
\label{int}
\int_\Omega \int_{A(\omega)}|g_i(\omega,x)|\varphi(\omega)(dx)\mu(d\omega)<\infty
\quad\mbox{for all}\quad i= 1,...,n
\end{equation}
and for some $\varphi\in\Phi.$ Put 
$g(\omega,x):= (g_1(\omega,x),...,g_n(\omega,x))$ for $(\omega,x)\in Gr(A)$ and
$$I_\varphi(g)(\omega):= \int_{A(\omega)}g(\omega,x)\varphi(\omega)(dx)\quad\mbox{and}\quad I_f(g)(\omega):= 
 g(\omega,f(\omega)), \ \omega\in \Omega,
$$
where $\varphi\in \Phi,$ $f\in F.$

Let $\cal G$ be a sub-$\sigma$-field of $\cal F.$  Following He and Sun \cite{hs} 
(see also Chapter XIV in  \cite{jacobs}), we say that 
$D\in \cal F$ is a $\cal G$-atom if $\mu(D)>0$ and for any $D_1\in \cal F$ there exists a set $D_2\in \cal G$ 
such that $\mu((D\cap D_1)\triangle (D\cap D_2) )=0.$\\

\noindent{\bf Remark 1} Inspired by the papers \cite{nh,r}, Dynkin and Evstigneev \cite{de}
used a different definition of a $\cal G$-atom in their studies of conditional expectation of correspondences 
(see page 334 in \cite{de}). However, as noted by He and Sun \cite{hs} (see page 54) the definitions of 
a $\cal G$-atom in  papers \cite{de,hs,jacobs} coincide.\\ 

We can now state our basic purification result.\\

\noindent{\bf Theorem 1} {\it Let the above assumptions on $(\Omega,\cal F,\mu)$,
the correspondence $A$ and the functions $g_i$ ($i=1,...,n)$ be satisfied. 
Let $\varphi\in \Phi$ be such that (\ref{int}) holds. 
If, in addition, the $\sigma$-field $\cal F$ has no $\cal G$-atoms, 
then there exists some $f\in F$ such that 
$$E(I_\varphi(g)|{\cal  G})=E(I_f(g)|{\cal G})\quad  (mod\ \mu). $$}

\noindent{\bf Remark 2} If ${\cal G}=\{\emptyset,\Omega\}$ and $A(\omega)=X$ for each $\omega\in \Omega,$
Theorem 1 reduces to a result closely related with Theorem 5.10 in Balder \cite{b4}, see also Theorem 2.1 in \cite{b5}. 
The proof given in \cite{b4}
makes use of the measurable correspondence $\Gamma$ from $\Omega$ to $X,$ where $\Gamma(\omega)$ 
is the intersection of all closed sets $C\subset X$ such that $\varphi(\omega)(C)=1.$ 
If we applied Balder's approach, then the sets $A(\omega),$ $\omega\in\Omega,$ would coincide with $\Gamma(\omega),$ 
which are in turn closed by his aasumption. In our case, 
 the sets $A(\omega)$ need not be closed. 
Therefore, the arguments given by Balder \cite{b4} 
must be modified appropriately to cover the case, where the sets  $A(\omega)$ depend on $\omega\in \Omega$ 
and may not be closed. The method used by Balder \cite{b4} is also connected with some measurability issues 
arising from an application of the measurable implicit function theorem.  Our proof
is somewhat simpler in the discussed special case of $\cal G.$ 
If  $\Omega$ is a Borel subset of a complete separable metric space and
${\cal G}=\{\emptyset,\Omega\}$, then our result is related with Theorem 1 
of Feinberg and Piunovskiy \cite{fp2}, see also Theorem
2.1 for one-step models in \cite{fp1}. The method used in \cite{fp1,fp2}
is based on Theorem 1 on ``strategic measures'' in \cite{fg} whose proof is long, quite involved and strongly based 
on the Borel structure of the spaces $\Omega$ and $X.$ Our approach is based on a relatively simple lemma 
stated below  and some standard measurable selection techniques. In the general case
with $\cal F$ having no $\cal G$-atom, we make use of the extension of
Lyapunov's theorem \cite{l} due to Dynkin and Evstigneev \cite{de}. 
 \\

Our proof is based on the following known fact.\\

\noindent{\bf Lemma 1} {\it Let $Y\in {\cal B}(X)$ and let $\nu\in P(X)$ be such that $\nu(Y)=1.$ 
Assume that $\psi:Y\to \mathbb{R}^n$ is a Borel measurable function such that $\int_Y\|\psi(x)\| \nu(dx)<\infty$ and 
define $R_\psi:=\{\psi(y): y\in Y\}.$ Then, $\int_Y\psi(x)\nu(dx)\in coR_\psi.$}\\

For a {\it proof} consider the distribution function of $\psi$ defined by: $p_\psi(B):=\nu(\psi^{-1}(B)), $ 
where $B \in {\cal B}(\mathbb{R}^n).$
Using Theorem 16.13 on page 229 in \cite{bill} and Lemma 3 on page 74 in \cite{f}, we obtain 
$$ \int_Y\psi(x)\nu(dx)=\int_{\mathbb{R}^n}yp_\psi(dy) \in coR_\psi.$$\\

\noindent{\bf Proof of Theorem 1} Let $\varphi\in \Phi$ be such that (\ref{int}) holds. For each $\omega\in \Omega,$ 
define $H(\omega):= \{g(\omega,y): y\in A(\omega)\}.$ By Lemma 1, we know that 
$I_\varphi(g)(\omega)\in coH(\omega)$ for all $\omega\in \Omega.$ Let $\Delta^{n+1}$ be the set of all probability vectors
in $\mathbb{R}^{n+1}.$ By Carath\'eodory's theorem, for any fixed $\omega\in\Omega$
there exist $(\lambda_1,...,\lambda_{n+1})\in \Delta^{n+1}$ and
$y_j \in A(\omega)$ ($j=1,...,n+1$) such that 
\begin{equation}
\label{K}
I_\varphi(g)(\omega)=\sum_{j=1}^{n+1}\lambda_j g(\omega,y_j).
\end{equation}
Let $X^{n+1}:= X\times\cdots\times X$ ($n+1$ times). Endow the metric space
$\Delta^{n+1}\times X^{n+1}$ with the Borel $\sigma$-field ${\cal B}(\Delta^{n+1}\times X^{n+1}).$
Let $K(\omega)$ be the set of all 
$(\lambda_1,...,\lambda_{n+1},y_1,...,y_{n+1})\in \Delta^{n+1}\times X^{n+1}$ for which (\ref{K}) holds. 
Using standard methods \cite{cv}, one can easily show that the graph of 
the correspondence $K$ belongs to the product $\sigma$-field 
${\cal F}\otimes {\cal B}(\Delta^{n+1}\times X^{n+1}).$ 
By the von Neumann-Aumann measurable selection theorem (see \cite{ab} or \cite{cv}), 
there exists measurable functions $\delta_j:\Omega \to [0,1]$ and $x_j:\Omega\to X$ ($j=1,...,n+1$) 
such that for all $\omega \in \Omega,$ we have 
$$(\delta_1(\omega),...,\delta_{n+1}(\omega), x_1(\omega),...,x_{n+1}(\omega))\in K(\omega).$$ 
Hence, we have
\begin{equation}
\label{e}
I_\varphi(g)(\omega)= \sum_{j=1}^{n+1}\delta_j(\omega)g(\omega,x_j(\omega))
\quad\mbox{for all}\quad \omega\in \Omega.
\end{equation} 
Put $h_j(\omega):=g(\omega,x_j(\omega)),$ $j=1,...,n+1,$ $\omega\in \Omega.$
Then, every $h_j$ is a measurable mapping from $\Omega$ into $\mathbb{R}^n$ and  (\ref{e}) can be rewritten as
\begin{equation}
\label{eq}
I_\varphi(g)(\omega)= \sum_{j=1}^{n+1}\delta_j(\omega)h_j(\omega)
\quad\mbox{for all}\quad \omega\in \Omega.
\end{equation} 
If $\cal F$ has no $\cal G$-atom, then by Remark 1 $\cal F$ has no 
$\cal G$-atom in the sense of Dynkin and Evstigneev \cite{de}. 
Therefore, we can apply the arguments given on pages 337-338 in \cite{de} to equation (\ref{eq}) 
and obtain a  partition of $\Omega$ into $n+1$ measurable subsets $\Gamma_1,...,\Gamma_{n+1}$ such that
\begin{equation}
\label{eqq}
E(I_\varphi(g)|{\cal G}) = E\left(\sum_{j=1}^{n+1} h_j1_{\Gamma_j}|{\cal G}\right)\quad (mod \ \mu).
\end{equation}
Here, $1_{\Gamma_j}$ is the indicator function of the set $\Gamma_j.$ 
Define $f(\omega):= x_j(\omega),$ whenever $\omega\in \Gamma_j.$ Then, $f\in F$ and from (\ref{eqq}), 
it follows that 
$$E(I_\varphi(g)|{\cal G}) =E(I_f(g)|{\cal G}) \quad (mod \ \mu),$$
which completes the proof.  $\Box$ \\

We close this section with some  corollary to Theorem 1 for a   
purification problem involving a finite family of measures. 
Although the approach is standard, it may have some application to  
the statistical decision theory \cite{b,dww1,dww2,fp2,f}. A similar result for Borel space 
$\Omega$ and possibly discontinuous integrands was formulated in Theorem 2 in \cite{fp2}.\\

\noindent{\bf Corollary 1} {\it Let $(\Omega, \cal F, \mu)$ be a complete probability space where 
$\mu= (\mu_1+\cdots +\mu_k)/k$ and every $\mu_i$ is a non-atomic probability measure ($i=1,...,k$). 
Assume that the graph $Gr(A)$ of $A$ belongs to the $\sigma$-field ${\cal F}\otimes {\cal B}(X).$
Let $\varphi\in \Phi$ and $\hat{g}(\omega,x) =(\hat{g}_1(\omega,x),...,\hat{g}_m(\omega,x))$ 
be a measurable mapping from $Gr(A)$ into $\mathbb{R}^m$ such that
$\int_\Omega \|\hat{g}(\omega,x)\|\varphi(\omega)(dx)\mu(d\omega)<\infty.$
Then, there exists some $f\in F$ such that
$$\int_\Omega I_\varphi(\hat{g})(\omega)\mu_i(d\omega)=
\int_\Omega I_f(\hat{g})(\omega)\mu_i(d\omega)\quad\mbox{for each}\quad i=1,...,k.$$}
\\

\noindent{\bf Proof} It is sufficient to apply Theorem 1 with 
${\cal G} =\{\emptyset,\Omega\}$ and the vector valued function
$g:Gr(A)\to \mathbb{R}^{km}$ given by 
$(\hat{g}\frac{d\mu_1}{d\mu},...,\hat{g}\frac{d\mu_k}{d\mu}).$
Here $\frac{d\mu_j}{d\mu}$ is the Radon-Nikodym derivative of $\mu_j$ with respect to $\mu.$ $\Box$

\section{An application to a Robust Variational Problem}

Let $c_i: Gr(A)\to \mathbb{R}_+:=[0,+\infty),$ $i=1,\ldots,m,$ and $u: Gr(A)\to \mathbb{R}$ be 
${\cal F}\otimes {\cal B}(X)$-measurable functions such that
$c_1(\omega,\cdot),\ldots,c_m(\omega,\cdot),$ $-u(\omega,\cdot)$ are 
inf-compact on $A(\omega)$ for each 
$\omega\in\Omega.$\footnote{Recall that inf-compactness of $c_i(\omega,\cdot)$ means that the set 
$\{x\in A(\omega):\ c_i(\omega,x)\le \ \beta\}$ is compact for each $\beta \in\mathbb{R}$.} 
In addition, suppose that there exists a  $\mu$-integrable function $\Lambda: \Omega\to \mathbb{R}_+$ such that
$$|u(\omega,x)|\le \Lambda (\omega)\quad \mbox{for}\quad (\omega,x)\in Gr(A).$$
Furthermore, we assume that the malevolent nature  chooses
a parameter $p$ from the given set  ${\cal P}.$ In this way, the nature has an influence on    
the distribution $q(\cdot|p)$ of a state of the world $\omega\in \Omega.$ Suppose further that 
$q(\cdot|p)$ is absolutely continuous with respect to $\mu$ for every $p\in{\cal P}.$ 
By  $\rho(\cdot,p)$ we denote the density of $q(\cdot|p).$ 

Let $a_1,\ldots,a_m$ be positive numbers. We study the following robust variational problem:
$$ (RVP)\qquad\qquad \sup_{f\in F}\inf_{p\in{\cal P}}\int_\Omega u(\omega,f(\omega))\rho(\omega,p)\mu(d\omega)$$
$$\mbox{subject to }\qquad \sup_{p\in{\cal P}} 
\int_\Omega c_i(\omega,f(\omega))\rho(\omega,p)\mu(d\omega)\le a_i, \quad i=1,\ldots,m.$$
By $\cal G$ we denote the sub-$\sigma$-field of $\cal F$ generated by the family of functions 
$\{\rho(\cdot,p)\}_{p\in{\cal P}}.$\\

\noindent{\bf Theorem 2} {\it Let $a_1,\ldots, a_m\in\mathbb{R}_+$ be such that the set of 
$\varphi\in\Phi$ for which
$$ \sup_{p\in{\cal P}}\int_\Omega\int_{A(\omega)} c_i(\omega,x)\varphi(\omega)(dx)\rho(\omega,p)\mu(d\omega)\le a_i,
\quad i=1,\ldots,m$$
is non-empty. Moreover, assume that the $\sigma$-field $\cal F$ has no $\cal G$-atoms.  
Then, there exists a non-randomised decision function $f_*\in F$ solving the problem (RVP).}\\

\noindent{\bf Proof } 
We partly follow the idea from \cite{b1}. By the Urysohn Metrisation Theorem (see Theorem 3.40 in \cite{ab}),
we embed $X$ in the Hilbert cube $H=[0,1]^\infty$ 
(equipped with the topology of pointwise convergence). Then, $X$ can be
identified with a dense Borel measurable subset of the compact metric space 
$\hat{X},$ where $\hat{X}$ is the closure of the
embedded space $X$ in $H.$ Observe that $Gr(A)$ is ${\cal F}\otimes {\cal B}(\hat{X})$-measurable. 
Next, we extend functions in the following way: $\hat{c}_i(\omega,x)=c_i(\omega,x)$ if $(\omega,x)\in Gr(A)$
and $\hat{c}_i(\omega,x)=+\infty$ if  $(\omega,x)\not\in Gr(A),$ $i=1,\ldots,m,$ and $\hat{u}(\omega,x)=u(\omega,x)$ 
if $(\omega,x)\in Gr(A)$
and $\hat{u}(\omega,x)=-\infty$ if  $(\omega,x)\not\in Gr(A).$ 
Note that  $\hat{c}_1,\ldots,\hat{c}_m,\hat{u}$ are ${\cal F}\otimes {\cal B}(\hat{X})$-measurable and
$\hat{c}_1(\omega,\cdot),\ldots,\hat{c}_m(\omega,\cdot),$ $-\hat{u}(\omega,\cdot)$ are 
inf-compact on $\hat{X}$ for every $\omega\in \Omega.$ Therefore, 
$\hat{c}_1(\omega,\cdot),\ldots,\hat{c}_m(\omega,\cdot),$ $-\hat{u}(\omega,\cdot)$ 
are lower semicontinuous on $\hat{X}$ for each $\omega\in\Omega.$ 

Let $\hat{\Phi}$ be the set of all measurable functions $\phi: \Omega\to P(\hat{X})$.  
Our original problem (RVP) we now replace 
by the convexified (relaxed) one:
$$ (CRVP)\qquad\qquad \sup_{\phi\in \hat{\Phi}}\inf_{p\in{\cal P}}\int_\Omega\int_{\hat{X}} 
\hat{u}(\omega,x)\phi(\omega)(dx)\rho(\omega,p)\mu(d\omega)$$
$$\mbox{subject to }\qquad \sup_{p\in{\cal P}} 
\int_\Omega \int_{\hat{X}} \hat{c}_i(\omega,x)\phi(\omega)(dx)\rho(\omega,p)\mu(d\omega)\le a_i, \quad i=1,\ldots,m.$$

Let $h:\Omega\times\hat{X}\to \mathbb{R}\cup\{+\infty\}$ 
be an ${\cal F}\otimes {\cal B}(\hat{X})$-measurable function.
The weak topology on $\hat{\Phi}$ is defined as the  
coarsest topology for which all functionals 
$$ \phi\to J^p_h(\phi):=\int_\Omega\int_{\hat{X}}h(\omega,x)\phi(\omega)(dx)\rho(\omega,p)\mu(d\omega)$$
are continuous for every integrably bounded Carath\'eodory function $h.$
Moreover, the weak topology on $\hat{\Phi}$ is also the coarsest topology for which the functionals $\phi\to J^p_h(\phi)$ 
are lower semicontinuous 
for every non-negative function $h$ such that $h(\omega,\cdot)$ is 
lower semicontinuous for every $\omega\in\Omega,$ see \cite{b1}. 
Hence, $\phi\to J^p_{\hat{c}_i}(\phi)$ is lower semicontinuous for every $i=1,\ldots,m$ and $p\in{\cal P}.$ 
Therefore, the function 
$$ \phi\to \sup_{p\in{\cal P}} J^p_{\hat{c}_i}(\phi)$$
is lower semicontinuous for each $i=1,\ldots,m.$
Additionally, $\hat{\Phi}$ is compact (see Theorem 1(i) in \cite{b2} or Theorem V.1 in \cite{cv}).
Consequently, the set
$$\hat{\Phi}_0:=\{\phi\in\hat{\Phi}:\  \sup_{p\in{\cal P}} 
J^p_{\hat{c}_1}(\phi)\le a_1,\ldots, \sup_{p\in{\cal P}} J^p_{\hat{c}_m}(\phi)\le a_m \}$$
is non-empty and compact. Put $\hat{u}_+=\max(\hat{u},0),$  $\hat{u}_-=\max(-\hat{u},0)$ 
and observe that $\hat{u}=\hat{u}_+-\hat{u}_- =u_+-\hat{u}_-.$
Therefore,
$$J^p_{\hat{u}}(\phi) = J^p_{u_+}(\phi)-J^p_{\hat{u}_-}(\phi) .$$
Since $u_+(\omega,\cdot)$ is upper semicontinuous on $\hat{X}$ and $\hat{u}_-(\omega,\cdot)$
is lower semicontinuous on $\hat{X},$ it follows from \cite{b1} that $\phi\to J^p_{u_+}(\phi)$ is upper semicontinuous and 
$\phi\to J^p_{\hat{u}_-}(\phi)$ is lower semicontinuous. 
Consequently, the mapping  $\phi\to J^p_{\hat{u}}(\phi)$ is upper semicontinuous for every 
$p\in{\cal P}$ and thus, the function
$$ \phi\to \inf_{p\in{\cal P}} J^p_{\hat{u}}(\phi)$$ is upper semicontinuous.
Hence,  there exists  $\phi_*\in \hat{\Phi}_0,$ which is a solution  to the problem (CRVP).

Now, due to our definition of extensions of functions used in (RVP),
we conclude that $\phi_*(\omega)(A(\omega))=1$ $\mu$-a.e. Let $N$  be the null set of all $\omega\in \Omega$ for which   
$\phi_*(\omega)(A(\omega))<1.$ Then, $N\in {\cal F}.$  Let $f$ be any element in $F.$ Define
$\varphi_*(\omega) := \phi_*(\omega)$ if  $\omega\in \Omega\setminus N$ and 
$\varphi_*(\omega)= \delta_{f(\omega)}$ if $\omega\in N,$
where $\delta_{f(\omega)}$ is a measure concentrated at the point $f(\omega).$ 
Thus, we infer that $\varphi_*\in \Phi$ and $\varphi_*$ is also a solution to (CRVP).

By Theorem 1, it follows that there exists, say 
$f_*\in F,$ such that
 $$E(I_{\varphi_*}(g)|{\cal  G})=E(I_{f_*}(g)|{\cal G})\quad  (mod\ \mu) ,$$
where $g_1:=u,$ $g_{k}:=c_{k-1}$ for $k=2,\ldots,n$ and $n:=m+1.$ Hence, for every $p\in{\cal P}$, we have
 $$E(I_{\varphi_*}(g)\rho(\cdot,p)|{\cal  G})=E(I_{\varphi_*}(g)|{\cal  G})\rho(\cdot,p) 
 =E(I_{f_*}(g)|{\cal G})\rho(\cdot,p)=
E(I_{f_*}(g)\rho(\cdot,p)|{\cal G})
\quad  (mod\ \mu).$$
Consequently, 
\begin{eqnarray*}
\int_\Omega I_{\varphi_*}(g)\rho(\omega,p)\mu(d\omega)&=&
\int_\Omega E(I_{\varphi_*}(g)\rho(\omega,p)|{\cal  G})\mu(d\omega)=\\
\int_\Omega E(I_{f_*}(g)\rho(\omega,p)|{\cal  G})\mu(d\omega)&=&\int_\Omega I_{f_*}(g)\rho(\omega,p)\mu(d\omega).
\end{eqnarray*}
Obviously, the value of the problem (RVP) is not greater than the value of (CRVP). 
On the other hand, by the above equalities we have
\begin{eqnarray*}
\lefteqn{ \sup_{\phi\in \hat{\Phi}}\inf_{p\in{\cal P}}\int_\Omega\int_{\hat{X}} 
\hat{u}(\omega,x)\phi(\omega)(dx)\rho(\omega,p)\mu(d\omega)}\\&&
=\inf_{p\in{\cal P}}\int_\Omega\int_{A(\omega)} 
u(\omega,x)\varphi_*(\omega)(dx)\rho(\omega,p)\mu(d\omega)\\
&&=\inf_{p\in{\cal P}}\int_\Omega
u(\omega,f_*(\omega))\rho(\omega,p)\mu(d\omega)
\le  \sup_{f\in F}\inf_{p\in{\cal P}} \int_\Omega
u(\omega,f(\omega))\rho(\omega,p)\mu(d\omega).
\end{eqnarray*}
Moreover, for every $i=1,\ldots,m$ we get
\begin{eqnarray*}
 a_i&\ge& \sup_{ p\in{\cal P}}\int_\Omega\int_{A(\omega)} 
c_i(\omega,x)\varphi_*(\omega)(dx)\rho(\omega,p)\mu(d\omega)\\
&=&\sup_{p\in{\cal P}}\int_\Omega
c_i(\omega,f_*(\omega))\rho(\omega,p)\mu(d\omega).
\end{eqnarray*}
This proves our assertion. $\Box$\\

Below we provide three examples of robust variational problems. In the first two examples we do not specify the functions 
$u$ and $c_1,\ldots, c_m.$\\

\noindent{\bf Example 1}  Define
$\Omega:=[0,1]\times [0,1]$ and assume that ${\cal F}$ be the completion of ${\cal B}([0,1])\otimes{\cal B}([0,1])$ 
with respect to the Lebesgue measure $\mu$ on the unit square.
Furthermore, assume that ${\cal P}=(0,1]$ and
$$\rho(\omega,p)=  
\left\{\begin{array}{l@{\quad\mbox{if}\quad}l}
 \frac 1p,& \omega_1\in[0,p], \ \omega_2\in[0,1]\\
0,&  \omega_1\in(p,1],\ \omega_2\in[0,1]
\end{array}\right.
\qquad\mbox{where}\quad    \omega=(\omega_1,\omega_2)\in\Omega,\ p\in {\cal P}.$$
Since $\rho(\cdot,p)$ is independent of $\omega_2,$ then $\sigma$-field generated by $\{\rho(\cdot,p)\}_{p\in{\cal P}}$ is
${\cal G}={\cal B}([0,1])\otimes \{\emptyset,[0,1]\}.$ Clearly, $\cal F$ has no $\cal G$-atom.\\

\noindent{\bf Example 2} 
Assume that $(\Omega,{\cal F},\mu)$ be a complete probability space with a non-atomic probability measure.
Let $\{B_j\}_{j\in\mathbb{N}}$
be a measurable partition of $\Omega,$ i.e., $\Omega=\bigcup_{j\in\mathbb{N}} B_j$ 
and $B_i\cap B_j=\emptyset$ for $i\not=j.$
 By $\cal P$ we denote the set of all density functions defined on this partition. Hence,
$p = \{x_j^p\}_{j\in\mathbb{N}} \in {\cal P}$ if 
$$\sum_{j=1}^\infty x_j^p\mu(B_j)=1\quad\mbox{and}\quad x_j^p\ge 0 \ \mbox{for every  }j\in\mathbb{N}.$$ 
Then, we define 
$$\rho(\omega,p):= \sum_{j=1}^\infty x_j^p\mbox{\bf 1}_{B_j}(\omega),$$
where $\mbox{\bf 1}_{B_j}$ denotes the characteristic function of $B_j\in {\cal F}.$
Since for each $p\in{\cal P}$  the function $\rho(\cdot,p)$ is constant on every set $B_j,$ $j \in \mathbb{N},$ it follows 
that the $\sigma$-field $\cal G$ generated by the family $\{\rho(\cdot,p)\}_{p\in{\cal P}}$ coincides 
with the  $\sigma$-field generated by the
partition  $\{B_j\}_{j\in\mathbb{N}}$. Clearly, $\cal F$ has no $\cal G$-atom.\\

Our last example is motivated by some issues in economic theory. Similar models with endogeneous shocks
were studied in dynamic stochastic games  \cite{duggan,hs}, where attention is only paid to randomized strategies 
of the players.\\

\noindent{\bf Example 3}  
Assume that  $\Omega:=Q\times R,$ where $Q$ and $R$ are complete separable metric spaces with their 
Borel $\sigma$-fields ${\cal B}(Q)$ and ${\cal B}( R).$ 
Consider a one-period model with a single firm. 
Let the state of an industry be $\omega=(k,r)\in\Omega,$ 
where $k\in Q$ is a capital stock and $r\in R$ is a random shock. Assume that the current state $\omega$ 
is chosen by the nature. Hence, ${\cal P}:=\Omega$ with its element $p=(k,r).$ 
The next state of the industry follows the distribution $q(\cdot|p)$ given by the formula
$$q(S|p)=\int_Q\int_R \mbox{\bf 1}_{S}(k',r')\nu(dr')\mu_Q(dk'|p), \quad S\in{\cal B}(Q)\otimes {\cal B}(R),$$
 where
\begin{itemize}
\item $\nu$ is a non-atomic probability measure on $(R,{\cal B}(R)),$
\item $\mu_Q(\cdot|p)$ is the marginal of $q(\cdot|p)$ on $Q;$ 
furthermore suppose that there exists a non-atomic probability measure $\lambda$ 
on $(Q,{\cal B}(Q))$ such that $\mu_Q(\cdot|p)$ is absolutely continuous with respect to 
$\lambda$ for all $p\in {\cal P};$ let $\rho(\cdot,p)$
be the corresponding Radon-Nikodym derivative.
\end{itemize}
Put $\mu:=\lambda\otimes \nu.$ 
Let ${\cal F}$ be the completion  of ${\cal B}(Q)\otimes{\cal B}(R)$ with respet to $\nu.$  
Then, the $\sigma$-field generated by the family
$\{\rho(\cdot,p)\}_{p\in{\cal P}}$ is ${\cal G}={\cal B}(Q)\otimes\{\emptyset,R\}.$
Since $\nu$ is non-atomic, $\cal F$ has no $\cal G$-atom under $\mu.$
 
The company has to choose a feasible production plan $a\in A(\omega')\subset \mathbb{R}^d$ in the next period, where
$\omega'=(k',r')\in \Omega.$ Assume that the function $u$ and $c_1,\ldots,c_m$ are ${\cal B}(Q)\otimes{\cal B}(R)
\otimes {\cal B}(\mathbb{R}^d)$ measurable. 
The objective of the firm is to maximize its future expected profit subject to the integral constraints corresponding to 
the expected production costs. 
In other words, the firm faces the following problem:
$$ (RVP)\qquad\qquad \sup_{f\in F}\inf_{p\in{\cal P}}\int_Q\int_R 
u((k',r'),f(k',r'))\rho((k',r'),p)\nu(dr')\lambda(dk')$$
$$\mbox{subject to }\qquad \sup_{p\in{\cal P}} 
\int_Q\int_R c_i((k',r'),f(k',r'))\rho((k',r'),p)\nu(dr')\lambda(dk')\le a_i, \quad i=1,\ldots,m.$$

\end{document}